\documentclass[12pt,thmsa]{article}
\usepackage{amsfonts}
\usepackage{amssymb}
\usepackage{amsthm}
\usepackage{latexsym}
\usepackage{amsmath}
\usepackage{amssymb}
\newtheorem{teo}{Theorem}[section]

\newtheorem{lema}[teo]{Lemma}
\newtheorem{cor}[teo]{Corollary}
\newtheorem{prop}[teo]{Proposition}

\newtheorem{obs2}[teo]{Remark}

\newtheorem{tea}{Theorem}[subsection]

\newtheorem{no2}[teo]{Note}

\newtheorem{no3}[tea]{Note}

\newenvironment{dem}{\begin{proof}[Proof]}{\end{proof}}

\newcommand{\Q}{\mathbb{Q}}
\newcommand{\disc}{{\rm disc}}

\newcommand{\F}{\mathbb{F}}

\newcommand{\PGL}{{\rm PGL}}

\newcommand{\GL}{{\rm GL}}
\newcommand{\N}{\mathbb{N}}

\title{Galois actions on Q-curves and Winding Quotients}
\author{Francesc Bars\thanks{supported by BFM2003-06092}\ \ and Luis  Dieulefait  \thanks{supported by
a MECD postdoctoral grant at the Centre de Recerca Matem{\'a}tica
from Ministerio de Educaci{\'o}n y Cultura}}

\begin{document}

\maketitle
\begin{abstract}
We prove two ``large images" results for the Galois
representations attached to a degree $d$ Q-curve $E$ over a
quadratic field $K$: if $K$ is arbitrary, we prove maximality of
the image for every prime $p >13$ not dividing $d$, provided that
$d$ is divisible by $q$ (but $d \neq q$) with $q=2$ or $3$ or $5$
or $7$ or $13$. If $K$ is real we prove maximality of the image
for every odd prime $p$ not dividing $d D$, where $D = \disc(K)$,
provided that $E$ is a semistable Q-curve. In both cases we make
the (standard) assumptions that $E$ does not have potentially good
reduction at all primes $p \nmid 6$ and that $d$ is square-free.
\end{abstract}


\section{Semistable Q-curves over real quadratic fields}

Let $K$ be a quadratic field, and let $E$ be a degree $d$ Q-curve defined over $K$.
Let $D= \disc(K)$. Assume that $E$ is semistable, i.e., that $E$ has good or
semistable reduction at every finite place $ \beta $ of $K$. Recall that we can
attach to $E$ a compatible family of Galois representations $\{ \sigma_\lambda \}$
of the  absolute Galois group of $\Q$: these representations can be seen as those
attached to the Weil restriction $A$ of $E$ to $\Q$, which is an abelian surface
with real multiplication by $ F:= \Q(\sqrt{ \pm d})$ (cf. [E]). Let us call $U$ the
set of primes dividing $D$. For
 primes not in $U$, it is clear that $A$ is also semistable, so in
particular for every prime $\lambda$ of $F$ dividing  a prime
$\ell$ not in $U$ the residual representation
$\bar{\sigma}_\lambda$ will be a representation ``semistable
outside $U$", i.e., it will be semistable (in the sense of [Ri
97]) at $\ell$ and locally at every prime $q \neq \ell$, $q
\not\in U$. This is equivalent to say that its Serre's weight will
be either $2$ or $\ell + 1$ and that the restriction to the
inertia groups $I_q$ will be unipotent, for every $q \neq \ell$,
$q \not\in U$ (cf. [Ri97]).\\
Imitating the argument of [Ri97], we want to show that in this
situation, if the image of $\bar{\sigma}_\ell$ is (irreducible
and) contained in the normalizer of a Cartan subgroup, then this
Cartan subgroup must correspond to the image of the Galois group
of $K$, i.e., the restriction to $K$ of $\bar{\sigma}_\ell$ must
be reducible. More precisely:

\begin{teo}
\label{teo:firstcase} Let $E$ be a semistable Q-curve over a
quadratic field $K$ as above. If $\ell \nmid 2 d D$, $\lambda \mid
\ell$, and $\bar{\sigma}_\lambda$ is irreducible with image
contained in the normalizer of a Cartan subgroup of $\GL(2,
\bar{\mathbb{F}}_\ell)$, then the restriction of this residual
representation to the Galois group of $K$ is reducible.
\end{teo}

\begin{dem} For any number field X, let us denote by $G_X$ its absolute
Galois group.\\
 We know that if we take $\ell \not\in U$ the
residual representation $\bar{\sigma}_\lambda$ is semistable
outside $U$. If this representation is irreducible and its image
is contained in the normalizer $N$ of a Cartan subgroup, then
there is a quadratic field $L$ such that the restriction of
$\bar{\sigma}_\lambda$ to $G_L$ is reducible and the quadratic
character $\psi$ corresponding to $L$
is a quotient of $\bar{\sigma}_\lambda$ (cf. [Ri 97]). \\
 Using the description of the restriction of
$\bar{\sigma}_\lambda$ to the inertia group $I_\ell$ in terms of
fundamental characters, and the fact that the restriction of
$\bar{\sigma}_\lambda$ to the inertia groups $I_q$, for every $q
\neq \ell$, $q \not \in U$, is unipotent, we conclude as in [Ri
97] that the quadratic character $\psi$ can only ramify at primes
in $U$, and therefore that the quadratic field $L$ is unramified
outside $U$, the ramification set of $K$. \\
On the other hand, we know (by Cebotarev) that the restriction to
$G_K$ of $\bar{\sigma}_\lambda$ is isomorphic to
$\bar{\sigma}_{E,\ell}$. Let us assume that
$\bar{\sigma}_{E,\ell}$ is irreducible (*). Its image is contained
in $N$, and since the restriction of $\bar{\sigma}_\lambda$ to
$G_L$ is reducible, it follows that the restriction of
$\bar{\sigma}_{E,\ell}$ to $G_{L\cdot K}$ is reducible. We are
again in the case of ``image contained in the normalizer of a
Cartan subgroup" but now for a representation of $G_K$. Once
again, the quadratic character $\psi'$ corresponding to the
extension $L \cdot K / K$ is a quotient of the residual
representation $\bar{\sigma}_{E,\ell}$. Using  the fact that the
curve $E$ is semistable we know that the restriction of this
residual representation to all inertia subgroups at places
relatively primes to $\ell$ give  unipotent groups, and this
implies as in [Ri97] that $\psi'$ is unramified outside (places
above) $\ell$.  But $\psi'$ corresponds to the extension $L \cdot
K / K$, and $L$ is unramified outside $U$, thus $\psi'$ is also
unramified outside (places above primes in) $U$. This two facts
entrain that $\ell \in U$, which is contrary to our hypothesis.\\
This proves that the assumption (*) contradicts the hypothesis of
the theorem, i.e., that the restriction to $G_K$ of
$\bar{\sigma}_\lambda$ is reducible, as we wanted.
\end{dem}

Keep the hypothesis of the theorem above, and assume furthermore that the field $K$
is real. Then, the conclusion of the theorem together with a standard trick show
that the image of $\bar{\sigma}_\lambda$ can not be (irreducible and) contained in
the normalizer of a non-split Cartan subgroup: the reason is simply that the
representation $\sigma_\lambda$ is odd, thus the image of $c$, the complex
conjugation, has eigenvalues $1$ and $-1$. In odd residual characteristic, this
gives an elements which is not contained in a non-split Cartan, but if we assume
that $K$ is real, we have $c$ contained in the group $G_K$, and we obtain a
contradiction because as a consequence of theorem \ref{teo:firstcase} the
restriction of $\bar{\sigma}_\lambda$ to $G_K$ must be contained in the  Cartan
subgroup. This, combined with Ellenberg's generalizations of the results of Mazur
and Momose (cf. [E]), shows that the image has to be large except for very
particular primes. In fact, we have the following:
 \begin{cor}
 \label{teo:corol}
  Let $E$ be a semistable Q-curve over a real
quadratic field $K$ of square-free degree $d$. Assume that $E$ does not have
potentially good reduction at all primes not dividing $6$. Then, if $D$ is the
discriminant of $K$, for every $\ell \nmid  d D$, $\ell
>13$ and $\lambda \mid \ell$, the image of the projective representation $P(\bar{\sigma}_\lambda)$
is the full $\PGL(2, {\mathbb{F}}_\ell)$.
 \end{cor}

\begin{section}{Q-curves of composite degree over quadratic fields}

Let $E$ be a $Q$-curve over a quadratic field $K$ of square-free
degree $d$. Let $\lambda$ be a prime of $K$ and let us consider
the projective representation $P(\overline{\sigma}_{\lambda})$
coming from $E$. We can characterize the image in a subgroup of
$PGL_2(\F_l)$ with $\lambda|l$ of the projective representation
$P(\overline{\sigma}_{\lambda})$ by points of modular curves as
follows (proposition 2.2 [E]):
\begin{enumerate}
\item $P(\overline{\sigma}_{\lambda})$ lies in a Borel subgroup,
then $E$ is a point of $X_0(dl)^K(\Q)$,
\item $P(\overline{\sigma}_{\lambda})$ lies in the normalizer of a
split Cartan subgroup then $E$ is a point of $X_0^s(d;l)^K(\Q)$,
\item $P(\overline{\sigma}_{\lambda})$ lies in the normalizer of a
non-split Cartan subgroup, then $E$ is a point of
$X_0^{ns}(d;l)^K(\Q)$;
\end{enumerate}
where $X^K(\Q)$ is the subset of $P\in X(K)$ such that
$P^{\sigma}=w_{d}P$ for $\sigma$ a generator of $Gal(K/\Q)$ where
$w_{d}$ is the Fricke or Atkin-Lehner involution.

We have the following results ([E], propositions 3.2, 3.4):
\begin{prop}\label{prop31} Let $E$ be a $Q$-curve of square-free degree $d$ over
$K$ a quadratic field. We have:
\begin{enumerate}
\item Suppose $P(\overline{\sigma}_{\lambda} )$ is reducible for some $l=11$ or
$l>13$ with $(p,d)=1$ where $\lambda|l$. Then $E$ has potentially
good reduction at all primes of $K$ of characteristic greater than
3. \item Suppose $P(\overline{\sigma}_{\lambda})$ lies in the
normalizer of a split Cartan subgroup of $PGL_2(\F_l)$ where
$\lambda|l$ for $l=11$ or $l>13$ with $(l,d)=1$. Then $E$ has good
reduction at all primes of $K$ not dividing 6.
\end{enumerate}
\end{prop}
After this result we need to study what happens when the image
lies in the non-split Cartan situation. For this case, Ellenberg
obtains for the situation of $K$ an imaginary quadratic field,
that there is a constant depending of the degree $d$ and the
quadratic imaginary field $K$ such that if the image of
$P(\overline{\sigma}_{\lambda})$ lies in a non-split Cartan and
$l>M_{d,K}$ then $E$ has potentially good reduction at all primes
of $K$, see proposition 3.6 [E]. He centers in the twisted version
for $X^K$ to obtain this result. We obtain a similar result in a
non-twisted situation for
$X^K$, and with $K$ non necessarily imaginary.\\

We impose once for all that $d$, the degree, is even. We denote
$d=2\tilde{d}$. First, let us construct an abelian variety
quotient of the Jacobian of $X_0^{ns}(2\tilde{d};l)$ on which
$w_{2\tilde{d}}$ acts as 1 and having $\Q$-rang zero. Then using
``standard" arguments, that we will reproduce here for reader's
convenience, we obtain our result on the non-split Cartan
situation.

By the Chen-Edixhoven theorem, we have an isogeny between
$J_0^{ns}(2;l)$ and $J_0(2l^2)/w_{l^2}$. Darmon and Merel [DM,
prop.7.1] construct an optimal quotient $A_f$ with $\Q$-rang zero.
They construct $A_f$ as the associated abelian variety to a form
$f\in S_2(\Gamma_0(2l^2))$ with $w_{l^2} f=f$ and $w_2f=-f$.

Then, in this situation, we construct now a quotient morphism
$$\pi_f:J_0(2\tilde{d}l^2)\rightarrow A_f'$$
such that  the actions of $w_{2\tilde{d}}$ and $w_{l^2}$ on
$J_0(2\tilde{d}l^2)$   give both the identity on $A_f'$ if
$\tilde{d}\neq 1$. Moreover, we can see that $A_f'$ is preserved
by the whole group $W$ of Atkin-Lehner involutions. We construct
$A_f'$ from $f\in S_2(\Gamma_0(2l^2) )$ and we go to increase the
level.

 We denote by $B_n$ the operator on modular forms of
weight 2 that acts as: $f|_{B_n}(\tau)=f(n\tau)=n^{-1}f|_{A_n}$,
where $A_n=\left(\begin{array}{cc}
n&0\\
0&1\\
\end{array}
\right)$ from level $M$ to level $Mk$ with $n|k$. We denote by
$$B_n:J_0(M)\rightarrow J_0(Mk)$$
the induced map on jacobians.

\begin{lema} With the above notation and supposing that $(\tilde{n},k)=1$
and $g$ is a modular form which is an eigenform for the
Atkin-Lehner involution $w_{\tilde{n}}$ in $J_0(M)$, then $g|_{B_{
{n}}}$ is also an eigenform for the Atkin-Lehner involution
$w_{\tilde{n}}$ in $J_0(Mk)$ with the same eigenvalue.
\end{lema}
\begin{dem} We only need to show that there exist $w_{\tilde{n},M}$ and $w_{\tilde{n},Mk}$, the Atkin-Lehner involution of
$\tilde{n}$ at level $M$ and $Mk$ respectively, such that:
$$A_n w_{\tilde{n},Mk}=w_{\tilde{n},M} A_n$$ which is easy to check.

\end{dem}

With the above lemma we can rewrite lemma 26 in [AL] as follows
\begin{lema}[Atkin-Lehner]\label{lemma26a} Let $g$ a form in $\Gamma_0(M)$, eigenform for all
the Atkin-Lehner involutions $w_l$ at this level. Let  $q$ be a
prime. Then the form
$$g|_{B_{q^{\alpha}}}\pm q^{(\delta-2\alpha)}g|_{B_1=Id}$$ is a form in
$\Gamma_0(Mq^{\alpha})$ which is an eigenform for all Atkin-Lehner
involutions at level $Mq^{\alpha}$ where $\delta$ is defined by
$q^{\gamma-\delta}||M$ and $q^{\gamma}||Mq^{\alpha}$. Moreover,
let us impose that $\delta\neq 2\alpha$. Then the eigenvalue of
this form for $w_{q^{v_q(Mq^{\alpha})}}$ is $\pm$ the eigenvalue
of $w_{q^{v_q(M)}}$ on $g$.
\end{lema}
\begin{obs2}[AL]\label{lemma26b} In the case $\delta=2\alpha$ let us take the form
$g|_{B_q^{\alpha}}$. Then it satisfies the following: it is an
eigenform for the Atkin-Lehner involutions at level $Mq^{\alpha}$
with eigenvalue for the Atkin-Lehner involution at $q$ equal to
that of the Atkin-Lehner involution at $q$ on $g$ ($g$ of level
$M$).
\end{obs2}
Let us remark that if  the condition $\delta\neq 2\alpha$ is
satisfied we can choose a form in level $Mq^{\alpha}$ with
eigenvalue of the Atkin-Lehner involution at $q$ as one wishes: 1
or -1. This condition is always satisfied if $(M,q)=1$, situation
that we will use in this article. With this remarks the following
lemma is clear by induction:
\begin{lema} Let $g$ be a modular form of level $M$ which is an
eigenvector for all the Atkin-Lehner involutions at level $M$.
Then we can construct by the above lemma a modular form $f$ of
level $Mk$ ($k\in\N$) which is an eigenvector for all the
Atkin-Lehner involutions at level $Mk$, and moreover the
eigenvalue at the primes $q|M$ with $(q,k)=1$ is the same that the
one for the Atkin-involution of this prime at $g$ at level $M$,
and we can choose (1 or -1) the eigenvalue for the Atkin-Lehner
involution at the primes $q$ with $(q,k)\neq1$ if this prime
satisfies the condition $\delta\neq 2\alpha$ of the above lemma.
\end{lema}
Let us write a result in the form that  will be usefull for our
exposition, noting here that the even level condition can be
removed.
\begin{cor} Let us write $\tilde{d}=p_1^{\alpha_1}\ldots p_r^{\alpha_r}$ with $(\tilde{d},2p^2)=1$. We have a map $$I_{\chi_{p_1},\ldots,\chi_{p_r}}:J_0(2p^2)\rightarrow
J_0(2\tilde{d}p^2)$$ whose image is stable under the action of
$W$, and we can choose the action of $w_{2\tilde{d}}$ on the
quotient as $\pm$ the action of $w_{2}$ for an initial form $g\in
S_2(\Gamma_0(2p^2))$ eigenform for the Atkin-Lehner involutions at
level $2p^2$.
\end{cor}
\begin{dem} From lemma 27 in [AL], we have a base for the modular forms
which are eigenforms for the Atkin-Lehner involutions. Applying
the lemma of Atkin-Lehner above we have the result for
$\tilde{d}=p_1^{\alpha_1}$, we have to consider
$I_{\chi_{p_1}}=|_{B_{p_1^{\alpha_1}}}+\chi(p_1)p_1^{-\alpha_1}|_{B_1=Id}$,
where we can choose $\chi(p_1)$  as 1 or -1
depending on how we want the Atkin-Lehner involution at the prime
$p_1$ to act on the quotient. Inductively we obtain the result.
\end{dem}

Applying the above corollary with $\tilde{d}$ square-free
($\alpha_i=1$) in our situation ($\tilde{d}\neq 1$) and choosing
$w_{2\tilde{d}}=1$, we take
$$A_f':=I_{\chi_{p_1},\ldots,\chi_{p_r}}(A_f),$$
which is by construction a subvariety of $J_0(2\tilde{d}l^2)$
isogenous to $A_f$ which is stable under $W$ (at level
$2\tilde{d}l^2$) on which $w_{2\tilde{d}}$ and $w_{l^2}$ acts as
identity. In particular the $\Q$-rank of $A_f'$ is zero (recall
that we started with
an $A_f$ of $\Q$-rank zero).\\

By the Chen-Edixhoven isomorphism, we obtain a quotient map
$$\pi_f':J_0^{ns}(2\tilde{d};l)\rightarrow A_f' .$$
$\pi_f'$ is compatible with the Hecke operators $T_n$ with
$(n,2\tilde{d}l)=1$ (see for example lemma 17 [AL]) and moreover
$\pi_f'\circ w_{2\tilde{d}}=\pi_f'$. Let us recall that we are
interested in points on $X_0^{ns}(2\tilde{d};l)^K(\Q)$ (we want to
study the non-split Cartan situation). We have the following
commutative diagram:
\begin{center}
$$\begin{array}{ccc}
J_0^{ns}(2\tilde{d};l)&\rightarrow&A_f'\\
\downarrow i&&\downarrow id\\
J:=J_0^{ns}(2\tilde{d};l)^K&\rightarrow&A_f'\\
\end{array}$$
\end{center}
where $i$ is an isomorphism such that
$i^{\sigma}=w_{2\tilde{d}}\circ i$ with $\sigma$ the non-trivial
element of $Gal(K/\Q)$. Observe that $\psi_{f}:=\pi_f'\circ
i^{-1}:J\rightarrow A_f'$ is defined over $\Q$ because,
$$\psi_f^{\sigma}=(\pi_f')^{\sigma}\circ(i^{-1})^{\sigma}=\pi_f'\circ
w_{2\tilde{d}}\circ i^{-1}=\pi_f'\circ i^{-1}=\psi_f.$$

Let $R_0$ be the ring of integers of $K(\zeta_l+\zeta_i^{-1})$ and
$R=R_0[1/2\tilde{d}l]$, then $X_0^{ns}(2\tilde{d};l)$ has a smooth
model over $R$ and the cusp $\infty$ of $X_0^{ns}(2\tilde{d};l)$
is defined over R [DM]. We define
$$h:X_0^{ns}(2\tilde{d};l)/R\rightarrow J_0^{ns}(2\tilde{d};l)/R$$ by
$h(P)=[P]-[\infty]$. Then it follows by lemma 3.8 [E]
\begin{lema} Let $\beta$ be a prime of $R$. Then the map,
$$\pi_f'\circ h:X_0^{ns}(2\tilde{d};l)/R\rightarrow A_f'/R$$ is a formal
immersion at the point $\overline{\infty}$ of
$X_0^{ns}(2\tilde{d};l)(\F_{\beta})$.
\end{lema}

 We can prove a result for the non-split Cartan situation with a
constant independent of the quadratic field.
\begin{prop}\label{prop38} Let $K$ be a quadratic field, and $E/K$ be
a $Q$-curve of square-free degree $d=2\tilde{d}$, with
$\tilde{d}>1$. Suppose that the image of
$P(\overline{\sigma}_{\lambda})$ lies in the normalizer of a
non-split Cartan subgroup of $PGL_2(\F_l)$ with $\lambda|l$ for
$l>13$ with $(2\tilde{d},l)=1$. Then $E$ has potentially good
reduction at all primes of $K$.
\end{prop}
\begin{dem} We can follow closely the proof of prop.3.6 in [E],
let us reproduce it here for reader's convenience. Take $\beta$ a
prime of $K$ where $E$ has potentially multiplicative reduction,
if $\beta| l$ then the image of the decomposition group
$G_{\beta}$ under $P(\overline{\sigma}_{\lambda})$ lies in a Borel
subgroup. By hypothesis this image lies in the normalizer of a
non-split Cartan subgroup. We conclude that the size of this image
has order at most 2, which means that $K_{\beta}$ contains
$\Q(\zeta_l+\zeta_l^{-1})$, which is impossible once $l\geq 7$.

Now let us suppose that $E$ has potentially multiplicative
reduction over $\beta$ with $\beta\nmid l$, denote by $l'$ the
prime of $\Q$ such that $\beta|l'$. It corresponds to a cusp on
$X_0^{ns}(2\tilde{d};l)$ where we will take reduction modulo
$\beta$. The cusps of $X_0^{ns}(2\tilde{d};l)$ have minimal field
of definition $\Q(\zeta_l+\zeta_l^{-1})$ [DM,\S5], and $K$ is
linearly disjoint from $\Q(\zeta_l+\zeta_l^{-1})$; it follows that
the cusps of $X_0^{ns}(2\tilde{d};l)$ which lie over $\infty\in
X_0(2\tilde{d})$ form a single orbit under $G_K$. If
$\tilde{\beta}$ is a prime of $L=K(\zeta_l+\zeta_l^{-1})$ over
$\beta$, then the point $P\in X_0^{ns}(2\tilde{d};l)(K)$
parametrizing $E$ reduces mod $\tilde{\beta}$ to some cusp $c$. By
applying Atkin-Lehner involutions at the primes dividing
$2\tilde{d}$, we can ensure that $P$ reduces to a cusp which lies
over $\infty$ in $X_0(2\tilde{d})$. By the transitivity of the
Galois action, we can choose $\tilde{\beta}$ so that $P$ actually
reduces to the cusp $\infty$ mod $\tilde{\beta}$. Using the
condition that a $K$-point of $X_0^{ns}(2\tilde{d};l)$ reduces to
$\infty$, we have then that the residue field
$\mathcal{O}_K/\beta$ contains $\zeta_l+\zeta_l^{-1}$, and this
implies that $(l')^4\equiv 1$ mod
$l$, in particular $l'\neq2,3$ when $l\geq 7$.\\
We have constructed a form $f$ and an abelian variety $A_f'$
isogenous to the one associated to $f$ with $\Q$-rank zero and
$w_{2\tilde{d}}$ acting as 1 on it, and we have a formal immersion
$\phi=\pi_f'\circ h$ at $\overline{\infty}$
$$X_0^{ns}(d;l)^K/R\rightarrow A_f'/R.$$
Let us consider $y=P$ our point from the $Q$-curve and $x=\infty$
at the curve $X=X_0^{ns}(2\tilde{d};l)/R_{\beta}$, we know that
they reduce at $\beta$ to the same cusp if $P$ has multiplicative
reduction. Let us consider then $\phi(P)$ the point in $A_f'(L)$
with $L=K(\zeta_l+\zeta_l^{-1})$. Let $n$ be an integer which
kills the subgroup of $J_0^{ns}(2\tilde{d};l)$ generated by cusps,
it exists by Drinfeld-Manin, then $n h(P)\in J_0^{ns}(2\tilde{d};
l)$ and let $\tau\in Gal(L/\Q)$ and not fixing $K$, then
$P^{\tau}=w_{2\tilde{d}}P$ and we obtain that
$n\phi(P)^{\tau}=n\phi(P)$ then lies in $A_f'(\Q)$ which is a
finite group and then torsion, concluding that $\phi(P)$ is
torsion (this is  getting a standard argument [DM, lemma 8.3]).

Since $l'>3$ the absolute ramification index of $R_{\beta}$ at
$l'$ is at most 2. Then it follows from known properties of
integer models (see for example [E, prop.3.1]) that $x$ and $y$
reduce to distinct point of $X$ at $\beta$, in contradiction with
our hypothesis on $E$.

\end{dem}
Putting together propositions \ref{prop31} and \ref{prop38}, we
obtain:
\begin{cor}\label{teo:corol2} Let $E$ be a Q-curve over a quadratic field $K$ of square-free
 composite degree $d=2\tilde{d}$, with $\tilde{d}>1$. Assume
 that $E$ does not have potentially
good reduction at all primes not dividing $6$. Then, for every
$\ell \nmid 2\tilde{d}$, $\ell
>13$ and $\lambda \mid \ell$, the image of the projective representation $P(\bar{\sigma}_\lambda)$
is the full $\PGL(2, {\mathbb{F}}_\ell)$.
\end{cor}

To conclude, observe that if we take a Q-curve over a quadratic
field whose degree $d$ is odd and composite (and square-free),
there are more cases where the above result still holds: for
example if $3 \mid d$ the result holds because all the required
results from [DM] (in particular, the existence of a non-trivial
Winding Quotient in $S_2(3 p^2)$) are also proved in this case.
Moreover, since the only property of the small primes $q=2$ or $3$
required for all the results we need from [DM] to hold is the fact
that the modular curve $X_0(q)$ has genus $0$, we can apply them
to any of $q= 2 , 3 , 5 , 7 , 13$, and so we conclude that the
above result applies whenever $d$ is composite (and square-free)
and divisible by one such prime $q$.

\end{section}

\section{Bibliography}\mbox{}\\

[AL] Atkin, A.O.L., Lehner, J., {\it Hecke operators on
$\Gamma_0(m)$}, Math. Ann., 185 (1970), 134-160.\\

[DM] Darmon, H., Merel, L., {\it Winding quotients and some
variants of Fermat's last theorem}, J. Reine Angew. Math., 490
(1997) 81-100.\\

[E] Ellenberg, J.,{\it Galois representations attached to Q-curves
and the generalized Fermat equation $A^4 + B^2 = C^p$}, preprint.\\

[ES] Ellenberg, J., Skinner, C., {\it On the modularity of Q-curves},
Duke Math. J., 109 (2001), 97-122.\\

[Ri97] Ribet, K., {\it Images of semistable Galois
representations}, Pacific J.  Math. {\bf 181} (1997), 277-297. \\

\vspace{2cm}

Francesc Bars Cortina, Depart. Matem\`atiques, Universitat
Aut\`onoma de Barcelona, 08193 Bellaterra. E-mail:
francesc@mat.uab.es \\

Luis Dieulefait, Depart. d'Algebra i Geometria, Facultat
Matem\`atiques, Universitat de Barcelona, Gran Via de les Corts
Catalanes 585, 08007 Barcelona. E-mail: luisd@mat.ub.es

\end{document}